%% file: C2.tex
\documentclass[11pt,a4paper]{article}

\usepackage[utf8]{inputenc}
\usepackage[T1]{fontenc}
\usepackage[francais]{babel}

\usepackage{concrete,euler}

\usepackage{amsmath}
\usepackage{amsfonts}
\usepackage{amssymb}
\usepackage{amsthm}
\usepackage{array}

\input{bemacros}

\input{bemathlatex}

\def\B{\mathcal{B}}
\def\L{\mathcal{L}}
\def\Ell{\mathcal{E}}

\def\O{\mathcal{O}}
\def\Sym{\mathfrak{S}}

\def\act{\cdot}
\def\inter{\cap}
\def\setdiff{\smallsetminus}
\def\h<#1,#2>{\langle #1,#2\rangle} 

\let\vol\abs

\def\surject{\twoheadrightarrow}
\def\inj{\hookrightarrow}
\def\einj{\mathrel{\inj\smash{\kern-1.1em\raisebox{.75ex}{$\scriptstyle \e$}}\kern0.5em}}
\def\identically{\equiv}

\author{Bernardo \textsc{Freitas Paulo da Costa}\footnote{%
Département de mathématiques de la faculté des sciences d’Orsay\hfil \break
Université Paris-Sud 11, 91405 Orsay CEDEX\hfil \break
\texttt{bernardo.da-costa@math.u-psud.fr}}}
\date{\today}
\title{Deux exemples sur la dimension moyenne d'un espace de courbes de Brody}

\begin{document}
\maketitle

\begin{abstract}
On étudie la dimension moyenne de l'espace de courbes $1$-Brody
à valeurs dans deux surfaces complexes : d'abord pour des surfaces de Hopf,
et ensuite pour $\P^2$ privé d'une droite.
On montre dans le premier cas que la dimension moyenne est nulle
via une borne sur la croissance des fonctions holomorphes
faisant apparaître le lemme de la dérivée logarithmique.
Pour montrer la positivité dans le deuxième exemple,
on relève de la droite à son complémentaire
un espace de courbes de Brody de dimension moyenne positive
construit par déformation d'une courbe elliptique.
\end{abstract}

\section*{Introduction}

La dimension moyenne a été développée par Gromov comme un invariant
capable de distinguer différents espaces de dimension topologique infinie~\cite{gromovmdim}.
Pour cela, on compare les effets de l'action d'un groupe $G$ sur ces espaces
avec la croissance du groupe.

On s'intéresse à l'espace des courbes entières
à valeurs dans une va\-ri\-é\-té compacte complexe $X$ et à dérivée uniformément bornée par 1.
Ces courbes sont dites \emph{de Brody} d'après le lemme de réparamétrisation dû à Brody~\cite{Brody}
qui produit de telles courbes entières.
On note cet espace $\B_1(X)$,
et on remarque qu'il admet une action naturelle de $\C$ par translation à la source.

Suivant l'approche métrique pour la dimension moyenne,
on munit $\B_1(X)$ de la distance $d(f,g) = \sup_{\abs{z} \leq 1} d_X(f(z),g(z))$,
qui induit la topologie de la convergence uniforme sur les compacts de $\C$.
L'action de $\C$ engendre (entre autres) la famille de distances
\[
d_R(f,g) = \sup_{\abs{z} \leq R + 1} d_X(f(z),g(z)),
\]
correspondant à la translation par tous les éléments d'un disque de rayon $R$.
La dimension moyenne de $\B_1(X)$ par rapport à cette action sera notée
$\beta(X)$ et vaut
\[
\lim_{\e \to 0} \lim_{R \to \infty} \frac{\dim_\e(\B_1(X),d_R)}{\pi R^2},
\]
où l'$\e$-dimension $\dim_\e(\B_1(X),d_R)$ est la plus petite dimension d'un polyèdre $P$
qui approche $\B_1(X)$ «~à $\e$ près~» :
il existe une application continue de $\B_1(X)$ vers $P$
dont les fibres sont de diamètre inférieur ou égal à $\e$ pour la distance $d_R$.

Si la définition de $\B_1(X)$ fait intervenir une métrique sur $TX$,
la compacité de $X$ entraîne que la nullité de $\beta(X)$
ne dépend pas de la métrique choisie.
De plus, on sait que $\beta(X) < \infty$,
et $\P^1$ est la seule courbe pour laquelle $\beta$ n'est pas nul.

\medskip

Dans cet article, on montre deux résultats, l'un de nullité, l'autre de positivité, pour $\beta(X)$.

Premièrement, on verra $W = \C^2 \setdiff \{(0,0)\}$ comme le revêtement
universel d'une surface de Hopf $\SH$,
muni de la métrique induite.
L'ordre de croissance, au sens de la théorie de Nevanlinna,
des courbes de Brody de $\SH$ relevées dans $W$ est alors inférieur ou égal à 1.
Par un raisonnement inspiré des arguments de~\cite[Lemma 2.1]{Tsukamoto-Moduli},
on montre que
\begin{thm}
La dimension moyenne de l'espace de courbes $1$-Brody à valeurs dans une surface de Hopf
est nulle.
\end{thm}
En particulier, pour les surfaces de Hopf fibrées sur $\P^1$ en courbes elliptiques,
(par exemple, $\SH = \raise.5ex\hbox{$W$} \big/ \raise-.5ex\hbox{$(z,w) \sim (2z, 2w)$}$
et $[(z,w)] \mapsto (z:w) \in \P^1$)
cela montre que la dimension moyenne peut devenir nulle par fibrations.

\medskip

Ensuite, on considère le cas des complémentaires
évoqué dans~\cite[Problem 1.14]{Tsukamoto-Moduli}.
Ce cas est intéressant parce que, même si ces complémentaires
ne contiennent pas de courbe rationnelle,
des courbes entières peuvent s'accumuler sur les diviseurs enlevés.

En outre, on possède une description assez précise des courbes de Brody dans
des complémentaires d'un «~grand~» nombre d'hyperplans en position générale
par des théorèmes à la Borel~\cite{Green-ThBorel}.
Ainsi, le complémentaire de 5 droites en position générale dans $\P^2$
est hyperbolique, n'admettant pas de courbe entière non-constante.
Pour le cas de 4 droites, les courbes entières sont
concentrées sur les diagonales reliant les points doubles.
En particulier, les courbes de Brody sont des exponentielles linéaires qui parcourent les $\C^*$ diagonaux.
C'est encore vrai que les courbes de Brody dans le complémentaire de 3 droites en position générale
dans $\P^2$ forment un espace de dimension finie,
puisque ses coordonnées homogènes sont encore des exponentielles linéaires~\cite{BertelootDuval}.

Les résultats récents d'Eremenko~\cite{eremenko-divisors},
généralisant le théorème de Clunie et Hayman~\cite{ClunieHayman},
montrent que l'ordre de croissance
d'une courbe de Brody à valeurs dans le complémentaire de 2 droites de $\P^2$ est au plus un.
Par~\cite[Theorem 1.9]{Tsukamoto-Moduli}, cela implique
que la dimension moyenne de l'espace de telles courbes est nulle,
même s'il s'agit maintenant d'un espace de dimension infinie.

On conclut la classification du point de vue de la dimension moyenne,
en montrant que
\begin{thm}
$\beta(\P^2 \setdiff \P^1) > 0$.
\end{thm}
\noindent où on définit $\B_1(\P^2 \setdiff \P^1)$
comme l'adhérence dans $\B_1(\P^2)$ des courbes $1$-Brody à valeurs dans $\P^2 \setdiff \P^1$.

Cela découle du théorème d'approximation suivant :
\begin{thm}
\label{thm:approx_elliptique}
On fixe une inclusion $\P^1 \inj \P^2$.
Soit $\psi : \C \to \P^1$ une courbe bipériodique de dérivée $\norm{d\psi} < 1$.
Alors il existe :
\begin{itemize}
\item un espace $M \subset \B_1(\P^1)$, de dimension moyenne positive,
      de fonctions uniformément proches de $\psi$ ;
\item pour toute fonction $\tilde\psi$ de $M$, une suite de courbes $1$-Brody
      à valeurs dans $\P^2 \setdiff \P^1$ convergeant, dans $\B_1(\P^2)$, vers $\tilde\psi$.
\end{itemize}
De plus, on peut choisir les courbes dans $\P^2 \setdiff \P^1 = \C^2$
parmi des relévées de $\tilde\psi$.
\end{thm}

Cet article répond à deux questions posées par Masaki Tsukamoto ;
je le remercie ici d'y avoir attiré mon attention.
Je remercie aussi Julien Duval, pour m'avoir transmis ces questions,
ses idées et nos discussions très fructueuses.

\section{Préliminaires}

\subsection{Courbes de Brody}

Une fonction holomorphe $f:\C \to X$ à valeurs dans une variété complexe
est dite une \emph{courbe de Brody} si la dérivée $df$
est uniformément bornée pour une métrique $\norm{\cdot}$ sur $X$.
On dit que $f$ est une courbe $c$-Brody lorsque $\norm{df} \leq c$.

Pour $X$ compacte, on note $\B_c(X)$ l'espace des courbes $c$-Brody à valeurs dans $X$,
muni de la topologie de la convergence uniforme sur les compacts, ce qui le rend compact.
Si $X = Y - D$ est le complémentaire d'un diviseur $D$ dans une variété compacte $Y$,
$\B_c(X)$ est défini comme l'adhérence dans $\B_c(Y)$ des courbes $c$-Brody à valeurs dans $X$.

La version suivante du lemme de Schwarz
assure que les petites déformations d'une courbe de Brody
restent encore de Brody :
\begin{lem}
\label{lem:schwarz}
Soit $X$ une variété complexe compacte.
Il existe $\e > 0$, dépendant uniquement de $X$,
tel que pour toute courbe $1$-Brody $f:\C \to X$,
les courbes $g:\C \to X$ $\e$-proches de $f$ sont $2$-Brody.
\end{lem}
\noindent On dit que $g$ est $\e$-proche de $f$ si, pour tout $z \in \C$,
$f(z)$ et $g(z)$ sont au plus à distance $\e$ sur $X$.
\begin{proof}
On se donne une collection finie de cartes locales de $X$,
donc tout ensemble de diamètre assez petit est inclus dans une carte locale.
L'image par $g$ des disques de centre $z$ et rayon $\delta$
est de diamètre borné par $\e + \norm{df}\delta$,
donc, pour $\e$ et $\delta$ assez petits, contenue dans une carte locale $\phi : U \to \C^n$.
On fixe un tel $\delta$.

Dans la carte $\phi(U)$, le lemme de Schwarz classique donne un contrôle en $\frac{\norm{d\phi} \cdot \e}{\delta}$
pour la différence entre les dérivées de $\phi \circ f$ et $\phi \circ g$ au centre du disque.
Revenant à $X$, on doit prendre en compte la distorsion des dérivées induite par la cartes,
aussi bien que le transport de la dérivée de $f(z)$ à $g(z)$.
Quitte à diminuer au départ les cartes locales, la distorsion est bornée ;
le transport est quant à lui donné par $(d\phi|_{g(z)})^{-1} \circ d\phi|_{f(z)}$,
qui est de norme proche de $1$ pour $g(z)$ proche de $f(z)$.
Comme $X$ est compact, toutes les estimées sont uniformes,
ce qui assure l'existence de $\e > 0$ suffisamment petit comme voulu.

Ce même argument montre qu'il existe $\e(\alpha)$ tel que l'on peut remplacer
$2$ par $1 + \alpha$ dans la conclusion.
\end{proof}

On rappelle dans la suite quelques notations et résultats de la théorie de distribution des valeurs sur $\P^1$.
Suivant Ahlfors, l'indicatrice de croissance de $f : \C \to \P^1$ est la fonction réelle
\[
T_f(r) = \int_0^r \frac{dt}{t} \int_{D_t} \norm{df}^2,
\]
où $D_t$ désigne le disque de rayon $t$ centré en $0$
et la norme de $df$ est calculée par la métrique de Fubini-Study.
Pour une fonction holomorphe $f : \C \to \C$,
\[
T_f(r) = \frac{1}{4\pi} \int_0^{2\pi} \log \left( 1 + \abs{f(r e^{i\theta})}^2 \right) \, d\theta.
\]
L'ordre de croissance d'une fonction $f$ est alors défini par
\[
\limsup_{r \to \infty} \frac{\log(T_f(r))}{\log r}.
\]
Ainsi, une courbe de Brody est d'ordre au plus $2$.

On note $n_f(r\,;a)$
le nombre de pré-images (avec multiplicité) de $a \in \P^1$ par $f$ dans le disque de rayon $r$,
et
\[
N_f(r\,;a) = \int_0^r \big( n_f(t\,;a) - n_f(0\,;a) \big) \frac{dt}{t} + n_f(0\,;a)\log r, 
\]
son intégrale logarithmique, appelée \emph{fonction d'impact}.
Écrivant $[a;b]$ pour la \emph{distance cordale} entre deux points de la sphère de Riemann,
\[
[a;b]^2 = \begin{cases} \frac{1}{(1 + \abs{a}^2)} & \text{si b = $\infty$} \\ \frac{\abs{a-b}^2}{(1 + \abs{a}^2)(1 + \abs{b}^2)} & \text{sinon} \end{cases}
\]
la \emph{fonction de proximité} est donnée par
\[
m_f(r\,;a) = \frac{1}{2\pi} \int_0^{2\pi} \log\left([f(re^{i\theta});a]^{-1}\right) \; d\theta.
\]
On a alors le premier théorème fondamental de Nevanlinna~\cite{Hayman} :
\[
N_f(r\,;a) + m_f(r\,;a) = T_f(r) + O(1).
\]
Comme $[a;b] \leq 1$, $m_f$ est positive, donc une fonction ne peut pas avoir asymptotiquement plus de zéros
que sa croissance le permet.

Par analogie avec les fonctions, on dit
qu'un ensemble discret $F$ du plan complexe est d'ordre
$\limsup_{r \to \infty} \frac{\log(N_F(r))}{\log r}$,
où $N_F(r)$ est défini comme $N_f$ avec $n_F(r) = \#(F \inter D_r)$ remplaçant $n_f$.
Comme $T_{f/g} \leq T_f + T_g + O(1)$, deux fonctions distinctes
ne peuvent pas coïncider sur un ensemble d'ordre strictement plus grand que leur ordre de croissance.

Enfin, le théorème de la dérivée logarithmique de Nevanlinna implique, en particulier,
que l'ordre de croissance de la dérivée $f'$
est au plus celui de la fonction de départ.

\subsection{Dimension moyenne}

On reprend les définitions de~\cite{gromovmdim} pour la dimension moyenne
dans le cas particulier des actions de $\C$.
\smallskip

Pour $(M,d)$ un espace métrique compact, son $\e$-dimension est
\[
\dim_\e(M,d) = \inf \big\{\, \dim P \mid P \mbox{ polyèdre tel qu'il existe } f: X \einj P \,\big\},
\]
où $f: X \einj P$ dénote une application continue $\e$-injective :
pour tout $p \in P$, le diamètre de $f^{-1}\{p\}$ est inférieur ou égal à $\e$.
La limite pour $\e \to 0^+$ de $\dim_\e(M,d)$ est la dimension topologique $\dim(M,d)$,
qui dépend uniquement de la topologie induite par $d$.

Si $\C$ agit sur $(M,d)$, on peut définir de nouvelles distances sur $X$,
d'abord la famille des translatées $d_t(x,y) = d(t\act x, t\act y)$ pour tout $t$ dans $\C$,
et ensuite les \emph{distances orbitales}
$d_\Omega = \sup_{t\in \Omega} d_t$
pour tout sous-ensemble borné $\Omega$ de $\C$.
Une suite $\Omega_n$ de sous-ensembles bornés est dite \emph{moyennable}
si l'aire des bords $w$-épaissis (où $w > 0$)
\[
\{\, t \in \C \mid \hbox{$t$ est à distance $\leq w$ de $\Omega_n$ et de son complémentaire} \,\}
\]
devient négligeable devant l'aire de $\Omega_n$ lorsque $n$ tend vers l'infini.
Un exemple est la suite $(D_n)$ des disques de rayon $n$ centrées en $0$.

La dimension moyenne de $(M,d)$ par rapport à une suite moyennable $\Omega_n$ est :
\[
\mdim(M,d : \{\Omega_n\}) = \lim_{\e \to 0} \lim_{n \to \infty} \frac{\dim_\e(M,d_{\Omega_n})}{\abs{\Omega_n}}.
\]
Quand l'action de $\C$ est uniformément continue, la limite intérieure existe
et ne dépend pas de la suite moyennable choisie,
et la limite en $\e \to 0$ ne dépend pas non plus de la distance $d$ de départ sur $M$.
On écrit alors simplement $\mdim(M:\C)$.

Comme pour la dimension topologique, la dimension moyenne est croissante.
C'est-à-dire, s'il existe une inclusion $\C$-équivariante $M \inj N$,
alors $\mdim(M:\C) \leq \mdim(N:\C)$.
De plus, si $P$ est (homéomorphe à) un polyèdre, $\mdim(P^\Z : \Z) = \dim(P)$.
Ce calcul et la propriété précédente sont au cœur des estimées de dimension moyenne.

\medskip

Sur $\C$ on possède un argument de plus, l'invariance d'échelle :
si $\Lambda = \Z z_1 \oplus \Z z_2$ est un réseau dans $\C$, alors
$\mdim(M:\Lambda) = \mdim(M:\C) \times \abs{\C / \Lambda}$.
La dimension moyenne par rapport à $\Lambda$ se construit de façon analogue :
on considère des sous-ensembles $\Omega_n \subset \Lambda$,
les distances orbitales induites,
et on remplace l'aire par la mesure de comptage.
De plus, comme $\Lambda$ est isomorphe à $\Z^2$,
on a encore $\mdim(P^\Lambda : \Lambda) = \dim(P)$ pour les polyèdres.

Cela permet de voir que la nullité de $\beta(X)$, la dimension moyenne de $\B_1(X)$,
est invariante par revêtements finis.
D'une part, si $X' \surject X$ est un revêtement, cela induit
le revêtement $\C$-équivariant $\B_1(X') \surject \B_ 1(X)$.
Par l'unicité après le choix d'un relèvement du point base,
on a $\mdim(\B_1(X'):\C) \leq \mdim(\B_1(X):\C)$
en choisissant $f : \B_1(X') \to X' : \phi \to \phi(0)$ dans le lemme suivant.
\begin{lem}
Soit $\pi : M' \to M$ continue et $\C$-équivariante entre deux espaces métriques compacts.
S'il existe un polyèdre $P$ et $f : M' \to P$ continue, commutant avec $\pi$ et telle que
$(\pi,f) : M' \to M \times P$ soit injective, alors
$\mdim(M':\C) \leq \mdim(M:\C)$.
\end{lem}
\begin{proof}
On considère une suite décroissante $\Lambda_n$ de sous-réseaux de $\C$,
et on construit les applications
\[
\pi_n: M' \to M \times (P)^{\Lambda_n} : m' \mapsto \big( \pi(m'), (f(g \act m'))_{g \in \Lambda_n} \big),
\]
qui sont $\Lambda_n$-équivariantes et injectives.
On en déduit $\mdim(M':\Lambda_n) \leq \mdim(M:\Lambda_n) + \dim(P)$,
et donc $\mdim(M':\C) \leq \mdim(M:\C) + \frac{\dim(P)}{\abs{\C/\Lambda_n}}$.
Faisant $n \to \infty$, on obtient l'inégalité voulue.
\end{proof}

D'autre part, on peut aussi contrôler la dimension du quotient par le
\begin{lem}
Soient $\Gamma$ un groupe fini dont l'action sur $M'$ commute avec celle de $\C$,
et $\pi : M' \map{/\raise -.5ex\hbox{$\scriptstyle\Gamma$}} M$ le quotient.
Alors
\[
\mdim(M : \C) \leq \#\Gamma \, \mdim(M' : \C).
\]
\end{lem}
\begin{proof}
Soient $\delta > 0$ et $\e > 0$ fixés.
Comme la limite en $\e$ dans la définition de la dimension moyenne est monotone,
pour $n$ assez grand il existe un polyèdre $P$
de dimension au plus $\vol{\Omega_n} (\mdim(M':\C) + \delta)$
et une fonction $\e$-injective $f: (M',d_n) \to P$.
On considère alors
\[
F: (M,d_n) \to P^k/\Sym_k : m \mapsto \left\{\, f(m') \in P \mid m' \in M' \text{ et } \pi(m') = m \,\right\} ,
\]
où $k = \#\Gamma$,
qui est $\e$-injective grâce à la compatibilité des métriques de $M'$ et $M$.
On a donc $\mdim(M:\C) \leq k \left(\mdim(M':\C) + \delta\right)$ pour tout $\delta$.
\end{proof}

\section{Surfaces de Hopf}

Une \emph{surface de Hopf} $\SH$ est une variété compacte complexe de dimension 2
dont le revêtement universel est $W = \C^2-\{(0,0)\}$.
Quand $\pi_1(\SH)$ est isomorphe à $\Z$, on montre~\cite{kodaira-ccasII}
que $\SH$ est difféomorphe à $S^1 \times S^3$ et
que la structure analytique peut être obtenue comme quotient de $W$
par l'action du groupe engendré par un difféomorphisme dilatant qui s'écrit~\cite{lattes-reduction2var}
\begin{equation*}
T_{(a,b,\lambda)} : W \to W : (z,w) \mapsto (az + \lambda w^n, bw),
\end{equation*}
où $a$, $b$ et $\lambda$ sont des nombres complexes, $|a| \geq |b| > 1$ ;
si $\lambda$ est différent de 0, on doit avoir $a = b^n$.
On les appelle \emph{surfaces de Hopf primaires}.
Les autres surfaces de Hopf, dites \emph{secondaires},
ont toujours une surface primaire comme revêtement fini.

Le lien établi par Tsukamoto entre la croissance des courbes de Brody
et la dimension moyenne~\cite[Theorem 1.9 et Lemma 2.1]{Tsukamoto-Moduli}
implique la remarque suivante :
si toutes les courbes de Brody dans une variété projective $X$
sont d'ordre strictement inférieur à $2$,
alors la dimension moyenne de $\B_1(X)$ est nulle.
Le raisonnement derrière la nullité de la dimension moyenne
sera de donner les détails et rendre précis l'argument suivant :

«~Les courbes de Brody à valeurs dans une surface de Hopf ont des coordonnées qui
sont des fonctions d'ordre au plus 1, donc forment un espace de dimension moyenne zéro.~»

Comme la nullité de la dimension moyenne ne change pas par revêtement fini,
on montrera le cas où $\SH$ est une surface de Hopf primaire.

\subsection{Analyse de croissance}
On prend sur $T\SH$ la métrique tirée en arrière de la métrique produit
par le difféomorphisme $\phi: \SH \to S^1 \times S^3$.
Par compacité, cela ne change pas l'analyse de $\beta(\SH)$.

Soit $p : W \to \SH$ la projection canonique.
Relevant $\phi \circ p : W \to S^1 \times S^3$ au revêtement universel,
on obtient une identification $\hat{\phi}$ entre $W$ et $\R \times S^3$,
équivariante par rapport à l'action du groupe fondamental,
qui rend commutatif le diagramme
\[
\begin{tabular}{r@{}ccc}
$\hat\phi :{}$& $W$   & $\to$ & $\R \times S^3$ \\
              & $\downarrow$ && $\downarrow$ \\
$\phi :{}$    &$\SH$  & $\to$ & $S^1 \times S^3$
\end{tabular}.
\]
On note $W_k$ la «~coquille sphérique~» correspondant dans $W$
au domaine fondamental $[k, k+1] \times S^3$.

Soit $f$ une courbe $1$-Brody sur $\SH$, de laquelle on fixe un relèvement $F=(f_1,f_2)$ à $W$.
Grâce au choix de la métrique de $\SH$, et vu que la projection $\pi_1 : S^1 \times S^3 \to S^1$
est à différentielle bornée par 1, on voit que
\[
\abs{d(\pi_1 \circ \phi \circ f)} \leq \abs{d\pi_1} \, \abs{d\phi} \, \abs{df} = \abs{df} \leq 1.
\]
Donc, si $\hat{\vphantom{\phi}\pi_1} \circ \hat\phi \circ F(0) \in [m_0,m_0+1]$,
alors l'image du disque centré en $0$ et de rayon $R$ par $\hat\phi \circ F$
est contenue dans $[m_0 - R, m_0 + 1 + R] \times S^3$.
De retour à $W$, cela veut dire que $F(D_R)$ est contenue dans
\[
\bigcup_{\abs{k} \leq R} T^k W_{m_0}.
\]
Soient alors $A_0$ et $B_0$ des bornes pour le module des coordonnées des points de $W_{m_0}$.
Traduisant l'inclusion précédente en coordonnées, on arrive à
\begin{align*}
|f_1| & \leq \abs{a}^R \left( A_0 + \abs{\lambda/a} R B_0^n \right) \\
|f_2| & \leq \abs{b}^R B_0 \qquad\qquad\qquad,
\end{align*}
car $T^R(z,w) = (a^R z + \lambda R a^{R-1} w^n, b^R w)$.
Cela montre que les coordonnées $f_i$ du relèvement sont des fonctions holomorphes
d'ordre au plus 1.

\subsection{Discrétisation}
Le quotient $f_2'/f_2$ est bien défini dans $\P^1$ dès que $f_2$ n'est pas identiquement nulle :
si $(\tilde{f}_1,\tilde{f}_2) = T^k(f_1, f_2)$ est un autre relèvement de $f$,
$\tilde{f}_2(z) = b^k f_2(z)$.
Cela revient à demander que $f$ ne soit pas incluse dans la courbe elliptique
$\Ell = \{ w = 0 \}$ de $\SH$.

Pour avoir une application bien définie sur $\B_1(\SH)$ entier,
nous introduisons $C(\P^1)$, le cône au-dessus de $\P^1$
réalisé par des rayons $r \cdot a$ pour $0 \leq r \leq \text{diam}(\SH)$ et $a \in \P^1$.
Cela permet de poser pour tout réseau $\Lambda \subset \C$ l'application de discrétisation
\[
\begin{array}{rccccccc}
P_\Lambda : & \B_1(\SH) &   \to   & \SH^\Lambda & \times & (C(\P^1))^\Lambda       \\
            &     f     & \mapsto & f |_\Lambda &    ,   & \left.\left( \text{dist}_\SH(f,\Ell) \cdot \frac{f_2'}{f_2} \right)\right|_\Lambda
\end{array}
\]
qui est continue et équivariante pour l'action de translation par $\Lambda$ à la source.

\medskip

On cherche à déterminer les réseaux $\Lambda$ qui donnent une discrétisation injective.
Soient $f$ et $g$ telles que $P_\Lambda(f) = P_\Lambda(g)$.
On a deux cas :
soit $f_2 \identically 0$,
soit $f_2'/f_2$ est bien définie.

Dans le premier cas, $g_2(\Lambda) \identically 0$,
et, comme c'est une application d'ordre 1, $g_2 \identically 0$.
Ainsi, les images de $f$ et $g$ sont contenues dans la courbe elliptique $\Ell$.
Les coordonnées $f_1$ et $g_1$ sont à valeurs dans $\C^*$ et d'ordre un, donc
$f_1(z) = \exp(Az + B)$ et $g_1(z) = \exp(Cz + D)$.
Comme $f = g$ sur $\Lambda$, pour un choix convenable
(et en fait presque tout choix) des générateurs de $\Lambda$
on obtient assez de contraintes sur les constantes $A$, $B$, $C$ et $D$
pour conclure que $f = g$ dans $\SH$.

Pour le second cas, comme $f_2$ est d'ordre 1,
il existe au plus un ensemble d'ordre 1 de points $z$ dans $\C$
où $f_2(z) = 0$, donc où $d(f(z),\Ell) = 0$.
Cela montre que $g_2(\Lambda)$ ne peut s'annuller que sur un ensemble d'ordre 1.
Pour tous les autres points de $\Lambda$, qui forment encore un ensemble d'ordre 2, on a
$f_2'/f_2 = g_2'/g_2$.
Mais ce sont deux applications d'ordre 1 (par le lemme de la dérivée logarithmique)
qui coïncident sur un ensemble d'ordre 2 ; elles sont donc égales sur $\C$ entier.
Intégrant cette égalité, on a $f_2 = C_2 g_2$ avec $C_2 \neq 0$.

Revenant à la projection sur $\SH^\Lambda$,
pour tout $z \in \Lambda$ il existe un entier $k(z)$ tel que
\[
(f_1,f_2)(z) = T^{k(z)}(g_1,g_2)(z).
\]
L'égalité en $0$ détermine $C_2 = b^{k(0)}$ car $f_2(0) = b^{k(0)} g_2(0)$.
La relation entre les deuxièmes coordonnées, à son tour, montre que $k(z) = k(0) = k$,
et alors $f_1(z) = a^k g_1(z) + k \lambda a^{k-1} g_2^n(z)$ pour tout $z$ de $\Lambda$.
Toutes les fonctions étant d'ordre 1, on a l'égalité précédente sur $\C$ entier,
donc $(f_1,f_2) = T^k(g_1,g_2)$, ce qui est exactement $f = g$.

\medskip

L'injectivité de $P_\Lambda$ montre que la dimension moyenne
de $\B_1(\SH)$ par rapport à l'action de $\Lambda$
est majorée par celle de $\SH^\Lambda \times (C(\P^1))^\Lambda$,
qui vaut $\dim(\SH) + \dim(C(\P^1)) = 7$.
Mais alors $\mdim(\B_1(\SH):\C) = \frac{\mdim(\B_1(\SH):\Lambda)}{|\C/\Lambda|}$,
et faisant tendre le covolume du réseau vers $+\infty$ on obtient
\[
\mdim(\B_1(\SH):\C) \leq \frac{7}{|\C/\Lambda|} \to 0.
\]

\section{Complémentaires}

Pour éviter de surcharger les notations, on se place dans le cas de $\P^2$
privé de la droite à l'infini, que l'on identifie à $\P^1$ de coordonnées $(x:y:0)$.
La stratégie est la suivante :
d'une part, on sait que la dimension moyenne de $\B(\P^1)$ est positive,
d'autre part, toute courbe dans $\P^1$ se relève à $\C^2$ qui est le complémentaire en étude.
On cherchera alors à montrer qu'il existe une partie $M \subset \B(\P^1)$,
de dimension moyenne positive, qui est dans $\B(\P^2 \setdiff \P^1)$.
Plus exactement, chaque courbe dans $M$ sera dans l'adhérence
d'une famille de ses relevées à $\C^2$, qui seront encore des courbes de Brody.

\subsection{Construction de $M$}
Cette partie exhibe $M$ comme un espace de déformations d'une courbe bipériodique dans $\P^1$.
Cela permettra de contrôler le relèvement des courbes vers $\C^2$
comme requis par le théorème~\ref{thm:approx_elliptique}.

\medskip

Soit $\psi : \C \to \P^1$ une application bipériodique satisfaisant $\norm{d\psi} < 1$.
Il existe alors une courbe elliptique $\Ell = \C / \Lambda$ et
d'une application $p : \Ell \to \P^1$ telles que $\psi = \pi \circ p$,
où $\pi: \C \to \Ell$ est la projection canonique.
On note $\L \to \Ell$ le fibré $p^* \O(1)$, et
$\sigma_0 = p^* X_0$ et $\sigma_1 = p^* X_1$ les sections des zéros et pôles de $p$.
Elles n'ont pas de zéro commun et $p$ s'écrit aussi $(\sigma_0 : \sigma_1)$.
De même, les $\psi_i = \pi^* \sigma_i$ sont des sections du fibré $\pi^* \! \L \to \C$,
et $\psi = (\psi_0 : \psi_1)$.

\medskip

Comme $\L$ est à courbure positive, il admet une métrique hermitienne
dont la forme de courbure est dans la classe de $\alpha dd^c \abs{z}^2$
pour une certaine constante $\alpha > 0$ ;
on peut même la choisir pour avoir exactement cette courbure.
Toute section $s$ trivialisant $\pi^*\! \L$ satisfait alors
$dd^c \log \norm{s}_\L = - \alpha dd^c \abs{z}^2$,
donc il existe une section holomorphe $f$ de norme
$\norm{f(z)}_\L = e^{-\alpha |z|^2}$.

Le fibré $\L$ donne des multiplicateurs
\[
e(\lambda, z) = \frac{\psi_0(z + \lambda)}{\psi_0(z)} \text{ pour $\lambda \in \Lambda$ et $z \in \C$,}
\]
correspondant à une translation par $\lambda$ sur la base de $\pi^* \! \L$
et au passage de la fibre $(\pi^* \! \L)_z$ vers $(\pi^* \! \L)_{z + \lambda}$
identifiées à $\L_{\pi(z)}$~\cite[p. 308]{GriffithsHarris}.
Cela permet de définir des nouvelles sections
$f_\lambda(z+\lambda) = e(\lambda, z) f(z)$ qui sont des translatées de $f$ ;
en particulier, $\norm{f_\lambda(z + \lambda)}_\L = \norm{f(z)}_\L$.

Ainsi, nous introduisons les convolutions de $f$ par une fonction bornée $a : \Lambda' \to \C$
\[
s_a(z) = \sum_{\lambda \in \Lambda'} a_\lambda \, f_\lambda(z)
\]
où $\Lambda'$ est un sous-réseau de $\Lambda$ qui sera déterminé à la fin de l'argument.
En dehors des zéros de $\psi_0$, on a l'expression
\[
\frac{s_a}{\psi_0}(z) = \sum_{\lambda \in \Lambda'} a_\lambda \, \frac{f}{\psi_0}(z - \lambda),
\]
qui peut être aussi utilisée pour définir $s_a$
et qui justifie le terme convolution.

\subsection{Contrôle des déformations}

On considère l'application suivante :
\[
\begin{array}{ccccc}
I: (D_\e)^{\Lambda'} &   \to   & H^0(\pi^* \! \L \times \pi^* \! \L)  &   \to   & \{\, f : \C \to \P^1 \text{ holomorphe} \,\} \\
           a         & \mapsto &        (\psi_0, \psi_1 + s_a)        & \mapsto &            (\psi_0 : \psi_1 + s_a)
\end{array}
\]
où on se restreint aux fonctions $a:\Lambda' \to \C$ bornées par $\e > 0$ à fixer.
Comme $\sigma_0$ et $\sigma_1$ n'ont pas de zéro commun,
$\norm{\psi_0}_\L^2 + \norm{\psi_1}_\L^2$ est minoré.
De plus, comme $\Lambda' \subset \Lambda$,
\[
\sup_z \norm{s_a(z)}_\L \leq \norm{a}_\infty \sup_z \left( \sum_{\lambda \in \Lambda} \norm{f(z+\lambda)}_\L \right) ,
\]
donc il existe $\e > 0$ et $R > 0$ tels que $(\psi_0, \psi_1 + s_a)$
évite un $R$-voisinage de la section nulle de $\pi^* \! \L \times \pi^* \! \L$
dès que $\norm{a} \leq \e$.
Cela assure que le passage à $\P^1$ est bien défini.

Comme la projection de $\C^2$ vers $\P^1$ est Lipschitzienne
dès que l'on évite un voisinage du centre de projection,
on voit que
\[
\psi_a = (\psi_0 : \psi_1 + s_a)
\]
est uniformément proche de $\psi = (\psi_0 : \psi_1)$,
donc, par le lemme de Schwarz, $\psi_a$ est 1-Brody pour $\e$ assez petit (car $\norm{d\psi} < 1$),
ce qui montre que l'image de $I$ est bien dans $\B_1(\P^1)$.

On relève alors les $\psi_a$ à $\C^2$ :
comme $f$ est une section jamais nulle, $\varphi_a = (\psi_0 : \psi_1 + s_a : f)$
ne coupe pas la droite à l'infini.
Par ailleurs, $M \varphi_a = (\psi_0 : \psi_1 + s_a : f/M)$ converge uniformément
vers $(\psi_0 : \psi_1 + s_a : 0)$ lorsque $M \to \infty$.
Donc, pour $M$ assez grand, les $M \varphi_a$ sont uniformément proches de
$\psi$, identifié à $(\psi_0 : \psi_1 : 0)$.
Et quitte à diminuer encore $\e$, les $M \varphi_a$ seront aussi $1$-Brody pour $M > M_0$,
ce qui conclut la deuxième partie du théorème~\ref{thm:approx_elliptique}.

\subsection{Dimension moyenne}

Il reste seulement à voir que $I : a \mapsto (\psi_0 : \psi_1 + s_a : 0)$
est injective, continue et $\Lambda'$-équivariante.
Cela montrera que son image $M$ satisfait
$\mdim(M:\Lambda') \geq \mdim\left( (D_\e)^{\Lambda'} : \Lambda' \right) = 2$,
et comme $M \subset \B_1(\P^2 \setdiff \P^1)$, on a la positivité annoncée.

La continuité découle du caractère linéaire de $a \mapsto s_a$
puis Lipschitzien de la projection de $\C^2 \setdiff D_R$ vers $\P^1$.
L'injectivité peut être obtenue par un choix convenable de $\Lambda'$ :
si $\Lambda'$ est tel que
\[
\sum_{0 \ne \lambda \in \Lambda'} \norm{f(\lambda)}_\L \leq \frac{1}{2} \norm{f(0)}_\L,
\]
alors la première application est injective puisqu'en choisissant $\ell$ tel que
$\abs{a_\ell} > (1 - \delta) \norm{a}_\infty$ on a
\begin{align*}
\norm{s_a(\ell)}_\L & \geq \abs{a_\ell} \, \norm{f(0)}_\L - \sum_{\ell \ne \lambda \in \Lambda'} \abs{a_\lambda} \, \norm{f(\ell - \lambda)}_\L \\
                    & > \left( \frac{1}{2} - \delta \right) \norm{a}_\infty \norm{f(0)}_\L.
\end{align*}
Mais $(\psi_0 : \psi_1 + s_a) = (\psi_0 : \psi_1 + s_b)$ si et seulement si
$s_a - s_b = 0$, car $\psi_0$ n'est pas identiquement nulle et toutes les fonctions sont holomorphes,
donc la composée $a \mapsto (\psi_0 : \psi_1 + s_a)$ est injective.

Enfin, la quasi-périodicité des sections de $\pi^*\! \L$
montre que $I : (D_\e)^{\Lambda'} \to \B_1(\P^1)$ est une inclusion $\Lambda'$-équivariante.
Pour toute translation $\tau_w$ par un élément $w \in \Lambda'$ on a
\begin{align*}
(\tau_w f_\lambda) (z) & = f_\lambda(z + w) = f_{\lambda-w}(z) e(w,z) \\
\noalign{et alors}
s_{(\tau_w a)}(z) & = \sum_{\lambda \in \Lambda'} a_{\lambda + w} \, f_\lambda(z) = \sum_{\lambda \in \Lambda'} a_\lambda \, f_{\lambda-w}(z) \\
                  & = \sum_{\lambda \in \Lambda'} a_\lambda \, f_\lambda(z+w) / e(w,z) = s_a(z+w) / e(w,z)
\end{align*}
et on conclut puisque $\psi_i(z+w) = \psi_i(z) e(w,z)$ pour $i = 0$ et~$1$.

\subsection{Remarques}
Pourvu que $k \geq 1$, on peut refaire cette construction dans $\P^n$ privé de $(n-k)$ hyperplans.
Leur intersection est un hyperplan de dimension au moins $k$,
donc une courbe elliptique $\Ell \to \P^k$ produit $k+1$ sections périodiques,
que l'on complète avec des sections bornées, jamais nulles et décroissantes.
Cela donne un espace $M^k \subset \B_1(\P^k)$ de dimension moyenne positive
de courbes de Brody qui sont dans le bord de $\B_1(\P^n \setdiff (n-k)\text{ hyperplans})$.

En outre, on peut ne pas passer à la limite des courbes dans $\C^2$ :
comme $\varphi = (\psi_0 : \psi_1 : f)$ est elle-même une courbe de Brody (disons $\norm{d\varphi} \leq c$),
le contrôle Lipschitzien des déformations dans $\P^2$ garantit que,
pour $\e$ suffisamment petit, il existe une image de $M$ au voisinage de $\varphi$,
donc cette image est dans $X_{c+1}$ par le «~lemme de Schwarz~».

On a vu que la dimension moyenne de l'espace de déformations ainsi construit
provient d'une partie de l'espace des courbes de Brody à valeurs dans
la droite enlevée.
Il est alors naturel de se poser les questions suivantes :
est-ce que toute courbe de Brody dans $\P^1$ est limite de courbes de Brody dans $\C^2$ ?
Est-ce que la dimension moyenne des courbes de Brody dans $\C^2$ est exactement
la dimension moyenne des courbes de Brody qui se relèvent à $\C^2$ ?

\bibliographystyle{alpha}
\NoAutoSpaceBeforeFDP
\bibliography{complexag,mdim,nevanlinna}

\end{document}

%% file: bemacros.tex
\def\map#1{\mathrel{\mathop{\longrightarrow} \limits^{#1}}}

\def\dim{\mathop{\rm dim}\nolimits}
\def\mdim{\mathop{\rm mdim}\nolimits}

\def\e{\varepsilon}


%% file: bemathlatex.tex
\def\Z{{\mathbf{Z}}}

\def\R{{\mathbf{R}}}
\def\C{{\mathbf{C}}}

\def\P{{\mathbf{P}}}

\def\SH{{\mathcal{S}}}

\def\del#1delz{\frac{\partial #1}{\partial z}}

\def\abs#1{{\lvert#1\rvert}}
\def\norm#1{{\left\lVert#1\right\rVert}}

\def\thname{Théorème}
\def\propname{Proposition}
\def\lemname{Lemme}
\def\faitname{Fait}
\newtheorem{thm}{\thname}

\newtheorem{lem}[thm]{\lemname}